\documentclass[10pt,twoside]{article}
\usepackage{Latex-document}
\almvol{00}
\almttone{Frontiers of Science Awards}
\almtttwo{for Math/TCIS/Phys}
\firstpage{1}
\usepackage[english]{babel}
\usepackage[utf8]{inputenc}

\usepackage{multicol}

\usepackage{tabu}

\usepackage{graphicx}
\usepackage{tikz}
\usetikzlibrary{patterns,shapes,decorations.pathmorphing,decorations.pathreplacing,calc,arrows}
\usepackage{verbatim}
\usepackage{amssymb}
\usepackage{amstext}
\usepackage{amsmath}
\usepackage{amscd}
\usepackage{amsthm}
\usepackage{amsfonts}
\usepackage{enumerate}
\usepackage{latexsym}
\usepackage{mathrsfs}
\usepackage{hyperref}
\usepackage{mathtools}
\usepackage[all]{xy}
\usepackage{pgfplots}
\usepackage{caption}
\xyoption{all}

\usepackage{pstricks}
\usepackage{lscape}
\usepackage{comment}

\numberwithin{equation}{section}

\newtheorem{thm}{Theorem}[section]
\newtheorem{defi}[thm]{Definition}

\newtheorem{cor}[thm]{Corollary}

\newtheorem{ques}[thm]{Question}
\newtheorem{rem}[thm]{Remark}

\DeclareMathOperator{\add}{\mathsf{add}}
\DeclareMathOperator{\moduleCategory}{\mathsf{mod}} \renewcommand{\mod}{\moduleCategory}\DeclareMathOperator{\per}{\mathsf{per}}
\DeclareMathOperator{\inj}{\mathsf{inj}}
\DeclareMathOperator{\proj}{\mathsf{proj}}
\DeclareMathOperator{\Mod}{\mathsf{Mod}}

\DeclareMathOperator{\Fac}{\mathsf{Fac}}

\newcommand{\TT}{\mathcal{T}}
\newcommand{\FF}{\mathcal{F}}
\newcommand{\Hom}{\operatorname{Hom}\nolimits}
\newcommand{\End}{\operatorname{End}\nolimits}
\newcommand{\Ext}{\operatorname{Ext}\nolimits}
\newcommand{\CC}{\mathcal{C}}
\newcommand{\DD}{\mathcal{D}}
\newcommand{\D}{\mathsf{D}}
\newcommand{\Db}{\mathsf{D}^{\rm b}}
\newcommand{\Kb}{\mathsf{K}^{\rm b}}
\newcommand{\R}{\mathbb{R}}
\newcommand{\N}{\mathbb{N}}
\newcommand{\op}{\operatorname{op}\nolimits}

\DeclareMathOperator{\Hasse}{\operatorname{H}\nolimits}
\DeclareMathOperator{\silt}{\mathsf{silt}}

\DeclareMathOperator{\sttilt}{\mathsf{s\tau\mbox{-}tilt}}
\DeclareMathOperator{\twosilt}{\mathsf{2\mbox{-}silt}}

\DeclareMathOperator{\ctilt}{\mathsf{c\mbox{-}tilt}}
\DeclareMathOperator{\tors}{\mathsf{tors}}
\DeclareMathOperator{\torf}{\mathsf{torf}}
\DeclareMathOperator{\ftors}{\mathsf{f\mbox{-}tors}}

\begin{document}

\markboth{\hfill{\rm Takahide Adachi, Osamu Iyama, Idun Reiten} \hfill}{\hfill {\rm On $\tau$-Tilting theory \hfill}}

\title{On $\tau$-tilting theory}

\author{Takahide Adachi, Osamu Iyama, Idun Reiten}

\begin{abstract}
We give a brief introduction to $\tau$-tilting theory \cite{AIR}. In particular, we will see how our theory unifies two different branches of tilting theory, namely, silting theory and cluster tilting theory. We also introduce the history and recent developments.
\end{abstract}

\maketitle

\setcounter{tocdepth}{1}
\tableofcontents

\section{History}

In this section, we give a quick review of tilting theory, including two predecessors of $\tau$-tilting theory, namely, silting theory and cluster tilting theory. We refer to \cite{AHK} for general introduction to various aspects of tilting theory. 

\subsection{Tilting theory}
Let $A$ be a ring. We denote by $\Mod A$ the category of $A$-modules, and by $\proj A$ the category of finitely generated projective $A$-modules.
We denote by $\D(A)$ the unbounded derived category of $\Mod A$, and $\per A$ the thick subcategory of $\D(A)$ generated by $A$. Then the canonical functor $\Kb(\proj A)\to\D(A)$ gives an equivalence $\Kb(\proj A)\simeq\per A$.
Two rings $A$ and $B$ are called \emph{derived equivalent} if $\D(A)$ and $\D(B)$ are equivalent as triangulated categories. This is equivalent to that $\per A$ and $\per B$ are equivalent as triangulated categories.

The following notion is basic.

\begin{defi}\label{define silting}
Let $\CC$ be a triangulated category. An object $T\in\CC$ is called \emph{tilting} if $\Hom_{\CC}(T,T[i])=0$ for all non-zero integers $i\neq0$ and $T$ generates $\CC$ as a thick subcategory. 
\end{defi}

For a ring $A$, a tilting object in $\per A$ is called a \emph{tilting complex} of $A$. The stalk complex $A\in\per A$ is a typical example of a tilting complex of $A$.
The derived equivalences between rings can be controlled by tilting complexes.

\begin{thm}\cite{R}
Two rings $A$ and $B$ are derived equivalent if and only if there exists a tilting complex $T$ of $A$ such that the ring $\End_{\D(A)}(T)$ is isomorphic to $B$.
\end{thm}

There is an important special class of tilting complexes. Recall that, for a collection $S$ of objects in an additive category, we denote by $\add S$ the full subcategory consisting of direct summands of finite direct sums of objects in $S$. A \emph{tilting $A$-module} \cite{BB,M} is an $A$-module $T$ such that $\Ext^i_A(T,T)=0$ holds for each $i\ge1$ and there exist exact sequences
\[0\to P_n\to\cdots\to P_0\to T\to0\ \mbox{ and }\ 0\to A\to T^0\to\cdots\to T^n\to0\]
with $n\ge0$, $P_i\in\proj A$ and $T^i\in\add T$ for each $i$. Tilting modules with $n=0$ are precisely the projective generators in Morita theory. Tilting modules with $n=1$ are often called \emph{classical tilting modules} and contain modules associated with BGP reflections for quiver representations \cite{BGP,APR}. 

We will review the development of tilting theory in the last 30 years, which is centered around the following question.

\begin{ques}
What kind of structures does the set of tilting complexes of $A$ have?
\end{ques}

\emph{Mutation} is a categorical operation that replaces a direct summand of a given tilting module/complex to obtain a new one. 
It was initiated by Riedtmann-Schofield and Happel-Unger \cite{RS,HU}, and also appeared independently in modular representation theory of finite groups \cite{O}.
However, mutation is often impossible inside the class of tilting modules/complexes. Silting theory and cluster tilting theory explained below address this important shortcoming.

\subsection{Silting theory}

The notion of silting objects is a natural generalization of tilting objects. Mutation is always possible in the larger class of silting objects. Moreover silting objects parametrize other important structures of the derived categories (see Section \ref{section 3.1}). This notion at least goes back to the study of t-structures of the derived categories in \cite{KV} and \cite{HKM}. In the rest of this subsection, we recall some basic results in \cite{AI}.

\begin{defi}\label{define silting}
Let $\CC$ be a triangulated category. 
An object $T\in\CC$ is called  \emph{presilting} if $\Hom_{\CC}(T,T[i])=0$ for all positive integers $i$. 
An object $T\in\CC$ is called \emph{silting} if it is presilting and $T$ generates $\CC$ as a thick subcategory.
\end{defi}

A typical example of a silting object is given by a connective differential graded algebra $A$, that is, the category $\per A$ has a silting object $A$.

An object $X$ in a Krull-Schmidt category is called \emph{basic} if it is a direct sum of pairwise non-isomorphic indecomposable objects. We denote by $|X|$ the number of non-isomorphic indecomposable direct summands of $X$.
Assume that $\CC$ is a Krull-Schmidt triangulated category.
We denote by $\silt\CC$ the set of isomorphism classes of basic silting objects of $\CC$. 
For $T,U\in\silt\CC$, we write $T\ge U$ if $\Hom_{\CC}(T,U[i])=0$ holds for all positive integers $i$. Then we obtain a partially ordered set $(\silt\CC,\ge)$. We denote by $\Hasse(\silt\CC)$ the Hasse quiver of $\silt\CC$. Thus the set of vertices is $\silt\CC$, and we draw an arrow $T\to U$ if $T>U$ and there does not exist $V\in\silt\CC$ satisfying $T>V>U$.

We recall the following general notion from, for example, \cite{AR}. Let $\CC$ be an additive category and $\DD$ a full subcategory of $\CC$ satisfying $\DD=\add\DD$.
For $C\in\CC$, a morphism $f:D\to C$ is called a \emph{right $\DD$-approximation} of $C$ if $D\in\DD$ holds and for each morphism $g:D'\to C$ with $D'\in\DD$, there exists $h:D'\to D$ satisfying $g=fh$. It is called \emph{minimal} if each morphism $g:D\to D$ satisfying $fg=f$ is an automorphism.
Dually a \emph{(minimal) left $\DD$-approximation} is defined. We call $\DD$ \emph{functorially finite} if each object in $\CC$ admits a right $\DD$-approximation and a left $\DD$-approximation.

Now assume that $\CC$ is a Krull-Schmidt triangulated category such that the endomorphism ring of each object is Noetherian. 
Take $T=\bigoplus_{i=1}^nT_i\in\silt\CC$ with indecomposable $T_i$. For each $1\le i\le n$, there exist triangles
\[T_i^*\to S_i\xrightarrow{f_i}T_i\to T_i^*[1]\ \mbox{ and }\ T_i\xrightarrow{g_i}S'_i\to {}^*T_i\to T_i[1],\] 
where $f_i$ is a minimal right $(\add T/T_i)$-approximation of $T_i$ and $g_i$ is a minimal left $(\add T/T_i)$-approximation of $T_i$.

\begin{thm}\cite{AI}\label{silting mutation}
Let $\CC$ be a Krull-Schmidt triangulated category and $T,U\in\silt\CC$.
\begin{enumerate}[\rm(a)]
\item The Grothendieck group $K_0(\CC)$ of $\CC$ has a basis given by the classes $[T_i]$ with $1\le i\le n$. Thus $|T|$ is constant for all $T\in\silt\CC$.
\item Both $\mu_i^+(T):=(T/T_i)\oplus T_i^*$ and $\mu_i^-(T):=(T/T_i)\oplus {}^*T_i$ belong to $\silt\CC$.
\item We have $\mu_i^-\circ\mu_i^+(T)\simeq T$, $\mu_i^+\circ\mu_i^-(T)\simeq T$ and $\mu_i^+(T)> T>\mu_i^-(T)$.
\item There exists an arrow $T\to U$ in $\Hasse(\silt\CC)$ if and only if $U=\mu_i^-(T)$ holds for some $1\le i\le n$ if and only if $T=\mu_i^+(U)$ holds for some $1\le i\le n$.
\item If $T>U$, then there exist $1\le i,j\le n$ such that $T>\mu_i^-(T)\ge U$ and $T\ge\mu_j^+(U)>U$.
\end{enumerate}
\end{thm}

The operations $\mu_i^+$ and $\mu_i^-$ are called \emph{(silting) mutation}. As a consequence of (a) and (c) above, we obtain the following result immediately.

\begin{cor}\label{Hasse 2n regular}
Let $n:=|A|$ and $T\in\silt A$. In $\Hasse(\silt A)$, precisely $n$ arrows start at $T$, and precisely $n$ arrows end at $T$. 
Thus the Hasse quiver $\Hasse(\silt A)$ is $2n$-regular for $n:=|A|$.
\end{cor}

We refer to Corollary \ref{cor: 2-term silting} for more results.

\subsection{Cluster tilting theory}

Cluster categories and their cluster tilting objects were introduced in \cite{BMRRT} to categorify cluster algebras of Fomin and Zelevinsky \cite{FZ}. There exists a bijection between (reachable) cluster tilting objects in a cluster category and clusters in the corresponding cluster algebra, see \cite{K}. Notice that a more general notion of cluster tilting subcategories in exact or triangulated categories plays an important role in higher dimensional Auslander-Reiten theory \cite{I}.

Let $k$ be a field and $\CC$ a $k$-linear Hom-finite triangulated category. We call $\CC$ \emph{2-Calabi-Yau} if there exists a functorial isomorphism
\[\Hom_{\CC}(X,Y)\simeq D\Hom_{\CC}(Y,X[2])\]
for each $X,Y\in\CC$, where $D$ is the $k$-dual. Basic examples are given by cluster categories \cite{BMRRT,Am,Re} and the stable categories of preprojective algebras of Dynkin type \cite{GLS}.

\begin{defi}
Let $\CC$ be a 2-Calabi-Yau triangulated category, and $T\in\CC$. We call $T$ \emph{rigid} if $\Hom_{\CC}(T,T[1])=0$. We call $T$ \emph{cluster tilting} if $\add T=\{X\in\CC\mid\Hom_{\CC}(T,X[1])=0\}$.
\end{defi}

We denote by $\ctilt\CC$ the set of isomorphism classes of basic cluster tilting objects of $\CC$.
For each $T=\bigoplus_{i=1}^nT_i\in\ctilt\CC$ with indecomposable $T_i$, there exist triangles
\[T_i^*\to U_i\xrightarrow{f_i}T_i\to T_i^*[1]\ \mbox{ and }\ T_i\xrightarrow{g_i}U_i\to {}^*T_i\to T_i[1],\] 
where $f_i$ is a minimal right $(\add T/T_i)$-approximation of $T_i$ and $g_i$ is a minimal left $(\add T/T_i)$-approximation of $T_i$. 

\begin{thm}\cite[5.1]{IY}\label{cluster tilting mutation}
Let $\CC$ be a 2-Calabi-Yau triangulated category, and $T\in\ctilt\CC$. Then
$T_i^*\simeq {}^*T_i$ holds, and $\mu_i(T):=(T/T_i)\oplus T_i^*\simeq(T/T_i)\oplus {}^*T_i$ belongs to $\ctilt\CC$.
\end{thm}

We refer to Corollary \ref{cor: cluster tilting} for more results. The operation $\mu_i$ is called \emph{(cluster tilting) mutation}. 
This is an analogue of silting mutation, but there are different features.
The set $\silt\CC$ has a partial order and the mutation has two directions, while the set $\ctilt\CC$ is more symmetric since two mutations coincide.

There is a similar notion of mutation for exceptional sequences \cite{Ru}. We refer to \cite{AI,BRT} for comparisons between them.

\section{$\tau$-tilting theory}

In this section, we review $\tau$-tilting theory developed in \cite{AIR}, see also introductory articles \cite{IR,T}.
Our main notion of support $\tau$-tilting modules is defined by using the Auslander-Reiten translation $\tau$. We show that support $\tau$-tilting modules parametrize important landmarks in the module and derived categories, and unify silting theory and cluster tilting theory.
We also explain important notions in $\tau$-tilting theory such as the canonical partial order, mutation and completion.

\subsection{Definition}

Let $A$ be a finite dimensional algebra over a field $k$, and $\mod A$ the category of finitely generated $A$-modules.
We have dualities
$D:=\Hom_k(-,k):\mod A\leftrightarrow\mod A^{\op}$ and 
$(-)^*:=\Hom_A(-,A):\proj A\leftrightarrow\proj A^{\op}$. They give an equivalence $\nu:=D(-)^*:\proj A\simeq\inj A$
called the \emph{Nakayama functor}, where $\inj A$ is the category of finitely generated injective $A$-modules. For $X\in\mod A$ with a minimal projective presentation
$P_1\xrightarrow{f}P_0\xrightarrow{}X\xrightarrow{}0$,
we define the \emph{Auslander-Reiten translation} $\tau X\in\mod A$ by the exact sequence
\[0\xrightarrow{}\tau X\xrightarrow{}\nu P_1\xrightarrow{\nu f}\nu P_0.\]
Then $\tau$ gives a bijection between the isomorphism classes of indecomposable non-projective $A$-modules and those of indecomposable non-injective $A$-modules. In fact, it gives an equivalence $\tau:\underline{\mod} A\simeq\overline{\mod} A$ between the stable and costable categories of $A$. Moreover, there exist functorial isomorphisms $\underline{\Hom}_A(\tau^{-1}Y,X)\simeq D\Ext^1_A(X,Y)\simeq\overline{\Hom}_A(Y,\tau X)$ called the \emph{Auslander-Reiten duality}.

\begin{defi}\cite{AIR}\label{define tau-tilting}
Let $A$ be a finite dimensional algebra over a field $k$ and $M\in\mod A$.
We call $M$ \emph{$\tau$-rigid} if $\Hom_A(M,\tau M)=0$. 
We call $M$ \emph{$\tau$-tilting} if it is $\tau$-rigid and satisfies $|M|=|A|$.
We call $M$ \emph{support $\tau$-tilting} if there exists an idempotent $e\in A$ such that $M$ is a $\tau$-tilting $(A/(e))$-module.
\end{defi}

We denote by $\sttilt A$ the set of isomorphism classes of basic support $\tau$-tilting $A$-modules.
For example, the Auslander-Reiten duality implies that each $\tau$-rigid $A$-module $X$ is rigid (that is, $\Ext^1_A(X,X)=0$), and that classical tilting $A$-modules are precisely faithful $\tau$-tilting $A$-modules.
When $A$ is hereditary, support $\tau$-tilting modules are precisely support tilting modules introduced in \cite{IT}.

\subsection{Unification}

The notion of \emph{t-structure} in the derived category enables us to recover the original abelian category \cite{BBD}. Since each tilting complex gives a derived equivalence, it gives rise to a t-structure. A \emph{torsion pair} is shadow of a t-structure in the module category. In particular, each classical tilting module gives rise to a torsion pair. 

\begin{defi}
Let $A$ be a finite dimensional algebra over a field $k$. A full subcategory $\TT$ of $\mod A$ is called a \emph{torsion class} (respectively, \emph{torsionfree class}) if it is closed under extensions and factor (respectively, sub-) modules.
\end{defi}

We denote by $\tors A$ (respectively, $\torf A$) the set of torsion (respectively, torsionfree) classes in $\mod A$. Then we have a bijection
\[\tors A\simeq\torf A\ \mbox{ given by }\ \TT\mapsto\TT^\perp:=\{X\in\mod A\mid\Hom_A(\TT,X)=0\},\] and the inverse map is given by  $\FF\mapsto{}^\perp\FF:=\{X\in\mod A\mid\Hom_A(X,\FF)=0\}$.
A pair $(\TT,\FF)\in\tors A\times\torf A$ is called a \emph{torsion pair} if $\FF=\TT^\perp$, or equivalently, $\TT={}^\perp\FF$. In this case, $\TT$ is functorially finite if and only if so is $\FF$. It is known that $\TT\in\tors A$ is functorially finite if and only if there exists $M\in\TT$ satisfying $\TT=\Fac M:=\{X\in\mod A\mid\mbox{there exists a surjection $M^{\oplus i}\to X$ for some $i\in\N$}\}$.
We denote by $\ftors A$ the set of functorially finite torsion classes in $\mod A$.

The following bijection is a starting point of our paper \cite{AIR}.

\begin{thm}\cite[2.7]{AIR}
Let $A$ be a finite dimensional algebra over a field $k$. Then there exists a bijection
\begin{equation}\label{sttilt ftors}
\sttilt A\simeq\ftors A\ \mbox{ given by }\ M\mapsto\Fac M.
\end{equation}
The inverse map is given by sending $\TT\in\ftors A$ to a basic projective generator of the exact category $\TT$.
\end{thm}

Using the bijection \eqref{sttilt ftors}, we define a partial order $\ge$ on $\sttilt A$:  For $T,U\in\sttilt A$, we write $T\ge U$ if $\Fac T\supset\Fac U$. By restricting the bijection \eqref{sttilt ftors}, we obtain a bijection between the isomorphism classes of basic classical tilting $A$-modules and $\{\TT\in\ftors A\mid DA\in\TT\}$ due to Hoshino and Smal\o\ 
\cite[VI.6.5]{ASS}. 

There is also a bijection between $\sttilt A$ and a subset of $\silt A:=\silt(\per A)$ defined as follows.

\begin{defi}\label{define 2-silting}
Let $A$ be a finite dimensional algebra over a field $k$. An object in $\per A$ is called \emph{2-term} if it is isomorphic to some $(P^i,d^i)\in\Kb(\proj A)$ satisfying $P^i=0$ for all $i\neq-1,0$. We denote by $\twosilt A$ the subset of $\silt A$ consisting of all 2-term complexes.
\end{defi}

For example, $A$ and $A[1]$ belongs to $\twosilt A$.
Note that the study of 2-term complexes goes back at least to \cite{Au}.
The importance of 2-term tilting/silting complexes was pointed out by the Tsukuba school \cite{HKM}, where a certain variation of \eqref{sttilt ftors} was given.

\begin{thm}\cite[3.2, 3.9]{AIR}
Let $A$ be a finite dimensional algebra over a field $k$.  Then there exists an isomorphism of posets
\begin{equation}\label{twosilt sttilt}
\twosilt A\simeq\sttilt A\ \mbox{ given by }\ P\mapsto H^0(P).
\end{equation}
\end{thm}

Recall that a finite dimensional $k$-algebra $A$ is called \emph{2-Calabi-Yau-tilted} if there exists a 2-Calabi-Yau triangulated category $\CC$ with cluster tilting object $T$ such that $A\simeq\End_{\CC}(T)$ as $k$-algebras. In view of the following bijection, cluster tilting theory can be regarded as $\tau$-tilting theory for 2-Calabi-Yau-tilted algebras.

\begin{thm}\cite[4.1]{AIR}
Let $\CC$ be a 2-Calabi-Yau triangulated category with cluster tilting object $T$, and $A:=\End_{\CC}(T)$. Then there exists a bijection
\begin{equation}\label{ctilt sttilt}
\ctilt\CC\simeq\sttilt A\ \mbox{ given by }\ U\mapsto\Hom_{\CC}(T,U).
\end{equation}
\end{thm}

Note that there is also a direct bijection $\ctilt\CC\simeq\twosilt A$ \cite[4.7]{AIR}.

\subsection{Completion and mutation}

To give a more detailed account to $\tau$-rigid modules, we introduce the following notion.

\begin{defi}
Let $A$ be a finite dimensional algebra over a field $k$, and $\underline{M}=(M,P)$ a pair of $M\in\mod A$ and $P\in\proj A$.
We call $\underline{M}$ \emph{$\tau$-rigid} if $M$ is $\tau$-rigid and $\Hom_A(P,M)=0$.
A $\tau$-rigid pair $\underline{M}$ is called \emph{$\tau$-tilting} if $|M|+|P|=|A|$.
\end{defi}

Each $M\in\sttilt A$ can be uniquely extended to a $\tau$-tilting pair $(M,eA)$, where $e$ is an idempotent appeared in Definition \ref{define tau-tilting}. Thus we can regard $\sttilt A$ as the set of isomorphism classes of basic $\tau$-tilting pairs of $A$, where two pairs $(M,P)$ and $(M',P')$ are called \emph{isomorphic} if $M\simeq M'$ and $P\simeq P'$, and a pair $(M,P)$ is called \emph{basic} if both $M$ and $P$ are basic.
Define the \emph{direct sum} of pairs by $(M,P)\oplus(M',P'):=(M\oplus M',P\oplus P')$. 
A \emph{completion} of a basic $\tau$-rigid pair $\underline{M}$ is a basic $\tau$-tilting pair which has $\underline{M}$ as a direct summand. A completion of $\underline{M}$ is \emph{maximum} (respectively, \emph{minimum}) if it is the maximum (respectively, minimum) of all completions of $\underline{M}$ with respect to the partial order $\ge$ on $\sttilt A$.

\begin{thm}\cite[2.10, 2.13, 2.17]{AIR}\label{completion}
Let $\underline{M}=(M,P)$ be a basic $\tau$-rigid pair.
\begin{enumerate}[\rm(a)]
\item ${}^\perp(\tau M)\cap P^\perp$ belongs to $\ftors A$, and the corresponding $\tau$-tilting pair via \eqref{sttilt ftors} gives the maximum completion of $\underline{M}$.
\item $\Fac M$ belongs to $\ftors A$, and the corresponding $\tau$-tilting pair via \eqref{sttilt ftors} gives the minimum completion of $\underline{M}$.
\item The following conditions are equivalent.
\begin{enumerate}[\rm(i)]
\item $\underline{M}$ is $\tau$-tilting.
\item The maximum and the minimum completions of $\underline{M}$ coincide.
\item If a pair $\underline{N}$ satisfies that $\underline{M}\oplus\underline{N}$ is $\tau$-rigid, then $\underline{N}\in\add\underline{M}$.
\item If $N\in\mod A$ and $Q\in\proj A$ satisfy $\Hom_A(M,\tau N)=0=\Hom_A(N,\tau M)$ and $\Hom_A(P\oplus Q,M\oplus N)=0$, then $(N,Q)\in\add\underline{M}$.
\end{enumerate}
\end{enumerate}
\end{thm}

Part (a) of Theorem \ref{completion} generalizes a famous result by Bongartz \cite[VI.2.4]{ASS}. 
Thus the maximum completion is often called \emph{Bongartz completion}. 

An \emph{almost complete $\tau$-tilting pair} is a $\tau$-rigid pair $(M,P)$ satisfying $|M|+|P|=|A|-1$.
In this case, we have the following results.

\begin{thm}\cite[2.18, 2.30]{AIR}\cite{Z}\label{mutation}
Let $\underline{M}=(M,P)$ be a basic almost complete $\tau$-tilting pair.
\begin{enumerate}[\rm(a)]
\item $\underline{M}$ has precisely two completions, that is, the maximum completion $\underline{M}\oplus\underline{X}$ with $\underline{X}=(X,0)$ and the minimum completion $\underline{M}\oplus\underline{Y}$ with $\underline{Y}=(Y,Q)$.
\item $X$ and $Y$ in (a) are related by an exact sequence $X\xrightarrow{g} M'\to Y\to0$ such that $g$ is a minimal left $(\add M)$-approximation of $X$.
\end{enumerate}
\end{thm}

We define \emph{mutation} in $\sttilt A$ by $\mu_{\underline{X}}^-(\underline{M}\oplus\underline{X}):=\underline{M}\oplus\underline{Y}$ and $\mu_{\underline{Y}}^+(\underline{M}\oplus\underline{Y}):=\underline{M}\oplus\underline{X}$ (as in Theorems \ref{silting mutation} and \ref{cluster tilting mutation}).
Part (b) allows us to calculate the mutation $\mu^-$ by using left approximations. Then the mutation $\mu^+$ can be calculated by using the canonical anti-isomorphism $\sttilt A\simeq\sttilt A^{\op}$ of posets and the mutation $\mu^-$ in $\sttilt A^{\op}$, see \cite[2.32]{AIR}.

Theorem \ref{mutation} has been generalized in \cite{J}: For a basic $\tau$-rigid pair $\underline{M}$ of $A$, there exists a finite dimensional algebra $B$, which can be explicitly constructed, such that there exists a bijection between $\sttilt B$ and the set of isomorphism classes of completions of $\underline{M}$. For further results in this direction, see \cite{AP}.

Now we consider the Hasse quiver $\Hasse(\sttilt A)$ of the partially ordered set $\sttilt A$. For a basic $\tau$-tilting pair $\underline{M}=\bigoplus_{i=1}^n\underline{M}_i$ with indecomposable $\underline{M}_i$, we use the notations $\mu^+_i(\underline{M})$ and $\mu^-_i(\underline{M})$.
The following result is immediate from Theorem \ref{silting mutation} and the bijection \eqref{twosilt sttilt}, but a proof in the context of $\tau$-tilting theory was given in \cite{AIR}.

\begin{thm}\cite[2.33, 2.35]{AIR}\label{Hasse sttilt}
Let $A$ be a finite dimensional algebra over a field $k$, and $\underline{M},\underline{N}\in\sttilt A$.
\begin{enumerate}[\rm(a)]
\item There exists an arrow $\underline{M}\to\underline{N}$ in $\Hasse(\sttilt A)$ if and only if $\underline{N}=\mu_i^-(\underline{M})$ holds for some $1\le i\le n$ if and only if $\underline{M}=\mu_i^+(\underline{N})$ holds for some $1\le i\le n$.
\item If $\underline{M}>\underline{N}$, then there exist $1\le i,j\le n$ such that $\underline{M}>\mu_i^-(\underline{M})\ge \underline{N}$ and $\underline{M}\ge\mu_j^+(\underline{N})>\underline{N}$.
\end{enumerate}
\end{thm}

By Theorems \ref{mutation} and \ref{Hasse sttilt}, we immediately obtain the following result, which is interesting since the larger Hasse quiver $\Hasse(\silt A)$ is $2n$-regular by Corollary \ref{Hasse 2n regular}.

\begin{cor}
The Hasse quiver $\Hasse(\sttilt A)$ is $n$-regular for $n:=|A|$.
\end{cor}

The following notion is a tilting theoretic generalization of representation-finiteness.

\begin{defi}\cite{DIJ}
A finite dimensional algebra $A$ over a field $k$ is called \emph{$\tau$-tilting finite} if $\sttilt A$ is a finite set.
\end{defi}

There are a number of characterizations of $\tau$-tilting finiteness \cite{DIJ}. As an immediate consequence of Theorem \ref{Hasse sttilt}, we obtain the following result, which can be regarded as a tilting theoretic analogue of the connectedness of Auslander-Reiten quivers of ring-indecomposable representation-finite algebras.

\begin{cor}\cite[2.38]{AIR}
If $\Hasse(\sttilt A)$ has a finite connected component $C$, then $C=\Hasse(\sttilt A)$ holds.
In particular, if $A$ is $\tau$-tilting finite, then $\Hasse(\sttilt A)$ is connected.
\end{cor}

Notice that $\Hasse(\sttilt A)$ is not necessarily connected if $A$ is $\tau$-tilting infinite. In fact, for each positive integer $\ell$, there exists a finite dimensional algebra $A$ such that $\Hasse(\sttilt A)$ has precisely $\ell$ connected components \cite{AHIKM2}, see also \cite{D,Te}.
We do not know if $\Hasse(\sttilt A)$ can have infinitely many connected components.

\subsection{Application}
Now we apply our results to silting theory and cluster tilting theory.
As in the case of $\tau$-rigid pairs, we use the following terminology: 
Let $\CC$ be a 2-Calabi-Yau triangulated category with cluster tilting object $P$. A \emph{completion} of a basic rigid object $T\in\CC$ is an element of $\ctilt\CC$ which has $T$ as a direct summand. An \emph{almost complete cluster tilting object} of $\CC$ is a rigid object $T$ satisfying $|T|=|P|-1$.
Using Theorems \ref{completion}, \ref{mutation} and the bijection \eqref{ctilt sttilt}, we obtain the following results.

\begin{cor}\label{cor: cluster tilting}
Let $\CC$ be a 2-Calabi-Yau triangulated category with cluster tilting object $P$.
\begin{enumerate}[\rm(a)]
\item \cite[2.6, 3.7]{ZZ}\cite[4.5(b)]{AIR} Let $T$ be a basic rigid object in $\CC$. Then $T$ has at least one completion. Moreover, the following conditions are equivalent.
\begin{enumerate}[\rm(i)]
\item $T$ is cluster tilting.
\item $|T|=|P|$.
\item $T$ has a unique completion.
\item If $U\in\CC$ satisfies that $T\oplus U$ is rigid, then $U\in\add T$.
\end{enumerate}
\item \cite[5.3]{IY}\cite[4.5(a)]{AIR} Let $T$ be a basic almost complete cluster tilting object in $\CC$. Then $T$ has precisely two completions. Moreover, they are mutation of each other in the sense of Theorem \ref{cluster tilting mutation}.
\end{enumerate}
\end{cor}

Let $A$ be a finite dimensional algebra over a field $k$. A \emph{completion} of a basic 2-term presilting complex $T$ of $A$ is an element of $\twosilt A$ which has $T$ as a direct summand.
A completion of $T$ is \emph{maximum} (respectively, minimum) if it is the maximum (respectively, minimum) of all completions of $T$ with respect to the partial order $\ge$ on $\silt A$.
An \emph{almost complete 2-term silting complex} of $A$ is a 2-term presilting complex $T$ satisfying $|T|=|A|-1$.
Using Theorems \ref{completion}, \ref{mutation} and the bijection \eqref{twosilt sttilt}, we obtain the following results. 

\begin{cor}\label{cor: 2-term silting}
Let $A$ be a finite dimensional algebra over a field $k$.
\begin{enumerate}[\rm(a)]
\item \cite[5.4]{DF}\cite[2.16]{Ai} Let $T$ be a basic 2-term presilting complex of $A$. Then $T$ has the maximum and the minimum completions. Moreover, the following conditions are equivalent.
\begin{enumerate}[\rm(i)]
\item $T$ is 2-term silting.
\item $|T|=|A|$.
\item The maximum and the minimum completions of $T$ coincide.
\item If $U\in\per A$ satisfies that $T\oplus U$ is 2-term presilting, then $U\in\add T$.
\end{enumerate}
\item \cite[5.7]{DF}\cite[3.8]{AIR} Let $T$ be a basic almost complete 2-term silting complex of $A$. Then $T$ has precisely two completions, that is, the maximum and the minimum completions. Moreover, they are mutation of each other in the sense of Theorem \ref{silting mutation}.
\item The Hasse quiver $\Hasse(\twosilt A)$ of $\twosilt A$ is $n$-regular for $n:=|A|$.
\end{enumerate}
\end{cor}

\begin{rem}
Corollary \ref{cor: 2-term silting} was partially shown in an independent earlier work \cite[Section 5]{DF}. However, even at the technical level, it lacks key ingredients of tilting theory such as partial order, mutation in Theorem \ref{silting mutation}, and connection with torsion pairs, support $\tau$-tilting modules and cluster tilting theory.

A standard direct proof of Corollary \ref{cor: 2-term silting} in the context of silting theory can be found in \cite[2.16]{Ai} and \cite[Section 2]{Ki}.
Then, thanks to well-established theory involving the bijection \eqref{twosilt sttilt}, it is possible to recover Theorems \ref{completion} and \ref{mutation} from Corollary \ref{cor: 2-term silting}.
\end{rem}

We end this section by remarks on completions of (not necessarily 2-term) presilting complexes.
A \emph{completion} of a basic presilting complex $T$ of $A$ is an element of $\silt A$ which has $T$ as a direct summand.
It was shown in \cite{Ka,LZ} that a basic presilting complex does not necessarily have a completion.
On the other hand, if a basic presilting complex $T$ with $|T|=|A|-1$ has a completion, then there are infinitely many completions of $T$ and all of them are related by iterated mutations \cite[3.9]{IYa}.

\section{Further results}

In the last decade, there has been great progress in $\tau$-tilting theory thanks to serious works by a number of researchers.
In this subsection, we try to review some of them. Unfortunately, due to page limitations, the list is far from complete.

\subsection{The correspondences in module/derived categories}\label{section 3.1}

For a finite dimensional algebra $A$ over a field $k$, there exist bijections between certain important objects/subcategories in $\mod A$, $\per A$ and the bounded derived category $\Db(A):=\Db(\mod A)$ of $\mod A$.
It is well-known that $\tors A$ corresponds bijectively with intermediate t-structures of $\Db(A)$ \cite{HRS}. Also it was shown more recently that there are bijections between the following sets \cite{KY}.
\begin{enumerate}[\rm(i)]
\item $\silt A$.
\item The set of bounded co-t-structures in $\per A$.
\item The set of isomorphism classes of simple-minded collections in $\Db(A)$.
\item The set of algebraic t-structures in $\Db(A)$.
\end{enumerate}
The restrictions of these bijections together with other works and our bijections \eqref{sttilt ftors} and \eqref{twosilt sttilt} can be summarized as follows, where we refer to \cite{As1} for details.

\begin{thm}\cite{BY,IT,MS,As1}\label{bijections}
There exist bijections between the following sets.
\begin{enumerate}[\rm(i)]
\item $\sttilt A$.
\item $\ftors A$.
\item $\twosilt A$.
\item The set of intermediate bounded co-t-structures in $\per A$.
\item The set of isomorphism classes of 2-term simple-minded collections in $\Db(A)$.
\item The set of intermediate algebraic t-structures in $\Db(A)$.
\item The set of isomorphism classes of left-finite wide subcategories of $\mod A$.
\item The set of isomorphism classes of left-finite semibricks of $A$.
\item The set of isomorphism classes of right-finite wide subcategories of $\mod A$.
\item The set of isomorphism classes of right-finite semibricks of $A$.
\end{enumerate}
\end{thm}

For further results in this direction, see \cite{ES,IJY,G}.

\subsection{Recent developments}

For each finite dimensional algebra $A$, one can define a non-singular fan in the real Grothendieck group $K_0(\proj A)_\R$, called the \emph{$g$-fan}, which consists of the cones associated with $\tau$-rigid pairs. 
The idea goes back at least to \cite{H}, and appeared in \cite{DF,DIJ}, then studied in \cite{AHIKM1,AHIKM2,BPPW} more systematically. For $g$-fans in cluster algebra setting, see \cite{BCDMTY,PPPP}.
Regarding $K_0(\proj A)_\R$ as the space of stability conditions appearing in the geometric invariant theory for representations of $A$, one can define a wall-chamber structure of $K_0(\proj A)_\R$ \cite{B,BKT,BST,As,F}.
The notion of the \emph{TF-equivalence} \cite{As} shows the relationship between $\tau$-tilting theory, stability conditions, and the canonical decompositions of presentation spaces \cite{DF}, see \cite{AsI}.
The study of signed ($\tau$-)exceptional sequences, the ($\tau$-)cluster morphism category and the picture group is also very active \cite{BM,BM2,HI,IgT}.

The structure of the sets in Theorem \ref{bijections} and also a larger set $\tors A$ than $\ftors A$ was investigated in \cite{BCZ,DIRRT,IRTT,Kas}. It is a very active subject to give an explicit description of the sets given in Theorem \ref{bijections} and $\tors A$ for important classes of algebras $A$, for example, (graded) gentle and Brauer graph algebras \cite{AAC,APS,CD,JSW,PPP}, (generalized) preprojective algebras \cite{Miz,KM,IRRT,FG}, algebras appearing in group theory \cite{EJR,KK}, Lie theory \cite{ALS,W} and algebraic geometry \cite{Aug2,HaW,HW,IW}, and many others \cite{Ad,Ao,IZ,MSk}.

The notion of $\tau$-tilting finiteness and related notions (e.g.\ silting discreteness, $g$-tameness) has been further investigated in a large number of works including \cite{AMY,AHMW,Aug,BPP,PYK}. In particular, brick/$\tau$-tilting versions of Brauer-Thrall conjectures and other related problems has been studied in several works including \cite{Mo,MP,P,ST,STV}.

Silting theory was studied for larger classes of rings \cite{Ki,Gn,IK}. Also a ``non-compact'' version of the notion of 2-term silting complexes has been systematically studied in a series of works including \cite{AMV}. This allows us to control a much wider family of torsion classes, and also sheds new light on the study of $\tau$-tilting finite algebras \cite{AMV2,S}.

\address{Faculty of Global and Science Studies, Yamaguchi University, 1677-1 Yoshida, Yamaguchi 753-8541, Japan.\\
\email{tadachi@yamaguchi-u.ac.jp}}

\address{Graduate School of Mathematical Sciences, University of Tokyo, 3-8-1 Komaba Meguro-ku Tokyo 153-8914, Japan.\\
\email{iyama@ms.u-tokyo.ac.jp}}

\end{document}